\documentstyle[11pt]{article}
\newcommand{\bc}{\begin{center}}
\newcommand{\br}{\begin{right}}
\newcommand{\ec}{\end{center}}
\newcommand{\be}{\begin{equation}}
\newcommand{\ee}{\end{equation}}

\newcommand{\pr}{\parallel}

\newcommand{\rar}{\rightarrow}

\newcommand{\lrar}{\longrightarrow}

\newcommand{\bmax}{\mbox{max} \,}

\newcommand{\meq}{\geq}

\newcommand{\al}{\alpha}
\newcommand{\bet}{\beta}
\newcommand{\del}{\delta}
\newcommand{\eps}{\varepsilon}

\newtheorem{thm}{Theorem}[section]

\newtheorem{rem}{Remark}[section]
\newtheorem{prob}{Problem}[section]
\newtheorem{defi}{Definition}[section]
\newtheorem{cor}{Corollary}[section]

\vspace{2.5cm}

\begin{document}

\begin{center}
{\bf CONSTRUCTION OF A STABILIZING CONTROL AND SOLUTION TO A
PROBLEM ABOUT THE CENTER AND FOCUS FOR DIFFERENTIAL SYSTEMS WITH A
POLYNOMIAL PART ON THE RIGHT SIDE}
\end{center}

\vspace{0.5cm}

\bc {\bf I.M.Proudnikov} \ec \vspace{0.5cm} \bc {\em Russia,
St. Petersburg University, pim\underline{ }10@hotmail.com} \ec

\vspace{0.3cm}

{\bf Abstract} Stationary differential systems with polynomial right
sides are considered. Necessary and sufficient
conditions are formulated when a given domain is a domain of asymptotic
stability
and the origin of coordinates is either focus or center. The problem of
construction of a stabilizing control in the form of polynomial is studied.

\bigskip


\vspace{1.5cm}

\section{\bf Introduction}

\vspace{1cm}

   In the $19^{th}$ century the great French mathematician Henri Poincar\'e formulated
the problem of finding stability conditions for differential systems
without calculation being a solution. At the end of $19^{th}$  and at the beginning
of $20^{th}$  centuries the great Russian mathematician A.M.Lyapunov developed the
mathematical stabilization theory for differential systems. For this purpose
he developed two methods. The first one is based on the
characteristic numbers of the fundamental matrix. The second one
is based on construction of special functions called Lyapunov's functions
which have properties similar to distance from a considered point to the origin
of coordinates (it is supposed that the zero solution is stationary). Lyapunov's
methods received world recognition.

   In this paper, differential systems with a right polynomial part $f(.)$
of       $$ x=(x_1,x_2,...,x_n) \in \Re^n $$ are considered i.e.
$$
              f(x)=(f_1(x),f_2(x),..., f_n(x))^{*},
$$
where
$$ f_p(x) = \sum_{l_1,l_2,...,l_n \in I_p}
a_{l_1,l_2,...,l_n}^{(p)} x_1^{l_1} x_2^{l_2}...x_n^{l_n},
$$
 \( p \in 1:n, l_1,l_2, ... l_n \) are non-negative integers and \(* \)
is the transposition sign, $a^{(p)}_{l_1,l_2,...,l_n}$ are
real-valued numbers, $I_p$ is the set of degrees of the polynomial
$f_p(x).$

  To analyze the system
$$
               \stackrel {\cdot}{x} = f(x)
$$
a method is suggested which is different from Lyapunov's methods and based on
a system transformation idea, so that we are able to say
something definite about stability.

   From the technical point of view the important problem is
finding a domain of asymptotic stability.
Conditions are given so that any domain including the origin of coordinates is one
of them. Necessary and sufficient conditions are found
when the origin of coordinates is either focus or center.

The case is considered when the coefficients $a_{l_1,l_2,...,l_n}^{(p)}$
depend on t. Sufficient conditions are given for the problems formulated
above.

   The systems considered are interesting because asymptotic
stability of the zero solution is equivalent to the following
statement, that any solution of our system starting from any
point of some small region of the origin of coordinates tended to be zero;
this is not correct in general case.

   Further, we solve the problem of
constructing a stabilizing control in any given region of the
origin of coordinates for the system
$$ \stackrel
{\cdot}{x}=f(x,u), $$
where \( x \in \Re^n \) is a phase vector ,
\( u \in R^r \) is a control, \( f(x,u)=(f_1(x,u),f_2(x,u),...,
f_n(x,u))^{*} \) is a vector polynomial of $x$ and $u$ i.e. $$
f_p(x,u) = \sum_{l_1,l_2,...,l_n,m_1,m_2,...,m_r \in I_p}
a_{l_1,l_2,...,l_n,m_1,m_2, ...,m_r}^{(p)}
x_1^{l_1}x_2^{l_2}...x_n^{l_n}u_1^{m_1}u_2^{m_2}... u_r^{m_r}, $$
 \( p \in 1:n,\,\,\,\,l_1,l_2,...,l_n, m_1,m_2, ..., m_r \) are non-negative
integers, $a^{(p)}_{l_1,l_2,...,l_n,m_1,m_2,...,m_r}$ are
real-valued numbers, $I_p$ is the set of degrees of the polynomial
$f_p(x,u).$ We assume that the zero vector \( 0 = (0,0,...,0) \in
\Re^n \) is a solution of the system i.e. \( f(0,0)=0. \)

   The theorem was proved that any region including the origin
of coordinates can be made a region of asymptotic stability if we
choose a suitable control $u(x).$ The control $u(x)$ is chosen as
a polynomial with a degree not higher than the degree of the
vector-function $f(x,u)$ as a function of $x$.

   The formulated problem about asymptotic stability of a
system with a right polynomial part is very important for practical
applications in physics and technics and can be found in the book by
V.I.Zubov "The lectures on the control theory", 1975. p.60.

\vspace{1.5cm}

\section {\bf Domains of asymptotic stability}

\vspace{1cm}

   Let us consider the differential system
\be
                \stackrel {\cdot}{x} = f(x), \qquad \qquad \qquad
\label{1control}
\ee
where \( x=(x_1,x_2,...,x_n) \in \Re^n \) and
$$
              f(x)=(f_1(x),f_2(x),..., f_n(x))^{*}.
$$
The vector-polynomials $f_p(.)$ and their coefficients $a^p_{l_1,l_2,...,l_n}$
satisfy the conditions in the Introduction.

We assume that $f(x) \neq 0$ for all $x \neq 0$ in
some neighborhood of the origin of coordinates $0.$

We will call the order $deg(f_p(.))$ of the polynomial $f_p(\cdot)$ the
maximal degree of the
polynomial $f_p(.)$ in the variables $x_j, j \in 1:n,$ in totality i.e.
$$
deg(f_p(x))=\max_{l_1,l_2,...,l_n \in I_p} (l_1+l_2+...+l_n).
$$

   So if for $n=2$ and $l_{ij} \neq 0$
$$ f_1(x)=l_{11}x_1^2+l_{12}x_2^2+l_{13}x_1x_2, \,\,\,\,\,
f_2(x)=l_{21}x_1^3+l_{22}x_2+l_{23}x_3, $$
then the order of
$f_1(x)$ is equal to two and the order of $f_2(x)$ is equal to
three.

   We will call the order of the function $f(.)$ the maximal degree of
the polynomials $f_p(x),p \in 1:n,$ regarding the variables $x_j, j \in
1:n,$ i.e.
$$
deg(f(x))=\max_{p \in 1:n} \, deg(f_p(x)).
$$

   Consequently, the order of the function $f(.)$ for the example written
above is equal to $\max (2,3) =3$.

   We assume that the definitions of stability and asymptotic
stability are known (see \cite{zub1}-\cite{kras1}).

   \begin{defi} \cite{zub1}, \cite{zub2} A domain $D, 0 \in int D$,
consisting of whole trajectories of the system (\ref{1control}), is called
a region of asymptotic stability  if the limit
\be
 \pr x(t,x_0,t_0) \pr \rar 0
\label{2control} \ee when $t \rar \infty$ is fulfilled for any
initial point $x_0 \in D$ and any solution $x(\cdot),x(t_0)=x_0,$
of the system (\ref{1control}). In this case we will call the
system (\ref{1control}) asymptotic stable in D.
\end{defi}

   The following problem is important for practical applications.

\begin{prob} It is known that zero solution of the system
\be
\stackrel{\cdot}{x} = Ax,
 \label{3control} \ee
$A[n \times n], x \in \Re^n$, is asymptotic stable if all
eigen-values of the matrix $A$ have negative real parts. In this
case the system (\ref{3control}) is asymptotic stable in $\Re^n$.
\label{1probcalopt}

   Consider now the differential system
\be
   \stackrel {\cdot}{x} = Ax+ \varphi (x)       \qquad \qquad \qquad
\label{4control}
\ee
where $\varphi (\cdot)$ is a vector polynomial with degree not less than
two, $\varphi(0)=0.$

   What are the conditions on the coefficients of the vector-
polynomial $\varphi(\cdot)$ from (\ref{4control}) under which the
system (\ref{4control}) is asymptotic stable in a given domain
$D,0 \in int D,$ i.e. for any solution of the system
(\ref{4control}) and any initial point $x_0 \in D$ the limit
(\ref{2control}) is true?
\end{prob}

   The considered circle of questions includes the problem about center
and focus. This problem is formulated in the following way.

\begin{prob} Assume that the origin of coordinates (0,0) for the
system (\ref{3control}), $n=2, \\ x \in \Re^2 $ is the  center.
The necessary and sufficient conditions for this are that all
eigen-values of the matrix A are imaginary. By adding a vector
polynomial $\varphi(\cdot)$ (system (\ref{4control})) the center
(0,0) can be a focus. It is needed to find conditions on the
coefficients of the polynomial $\varphi(\cdot)$ that the point
(0,0) was the center of the differential system (\ref{4control}).
\label{2probcalopt}
\end{prob}

   Let us rewrite the system (\ref{1control}) in equivalent form
\be
              \stackrel{\cdot}{x} = A(x) x.
\label{5control}
\ee
   The elements $a_{ij}(x)$ of a matrix $A(x) [n \times n]$ are  continuous
polynomial functions of x. Conversion (\ref{4control}) to (\ref{5control})
 is not unique. It can be done in an infinite number of ways.
In fact, $$ a_{ij}(x)= \sum_{l_1+l_2+...+l_n \leq deg(f)}
\al^{(i)}_{j,l_1,l_2,...,l_n}(x) x_1^{l_1}x_2^{l_2}...
x_j^{l_j-1}... x_n^{l_n} $$  if $l_j \meq 1.$ We have the
following correlation for the coefficients
\be
\sum_{j}\al^{(i)}_{j,l_1,l_2,...,l_n}(x)=
a_{l_1,l_2,...,l_n}^{(i)} \label{5acontrol} \ee for all $x$ from
some region $D, 0 \in int D.$

For instance, let the system (\ref{1control}) have the right part
$$
          f(x) = \left( \begin{array}{r}
                         x_1^2x_2+x_2+2x_1x_2^2 \\
                         -x_1+3x_1^2x_2-2x_1x_2^2
                  \end{array} \right) .
$$
Then
$$
          A(x) = \left( \begin{array}{rr}
                        \al_1x_1x_2+\al_2x_2^2 & 1+\bet_1x_1^2+\bet_2x_1x_2 \\
                        -1+\gamma_1x_1x_2+\gamma_2x_2^2 & \del_1x_1^2+\del_2x_1x_2
                  \end{array} \right)
$$
where the coefficients $\al_i, \bet_i, \gamma_i, \bet_i, i=1,2 $ such that
$$
               \al_1+\bet_1=1, \,\,\,\,\, \gamma_1+\del_1=3,
$$
$$
                \al_2+\bet_2=2, \,\,\,\,\, \gamma_2+\del_2=-2.
$$
   The system (\ref{4control}) for the given vector function $f(\cdot)$
can be rewritten in the following form
\be
               \stackrel{\cdot}{x}=A_0 x+C(x)x
\label{5bcontrol} \ee
 where $$
                C(x)=\left( \begin{array}{rr}
                            \al_1x_1x_2+\al_2x_2^2 & \bet_1x_1^2+\bet_2x_1x_2       \\
                            \gamma_1x_1x_2+\gamma_2x_2^2 & \del_1x_1^2+\del_2x_1x_2
                      \end{array} \right),
$$
$$
                A_0 = \left( \begin{array}{rr}
                           0 & 1 \\
                           -1 & 0
                     \end{array} \right).
$$
   Since the eigen-values of the matrix $A_0$ are $\lambda _{1,2}=\pm i$, the
point (0,0) is the center for the linearized system. But it is
question for the nonlinear system (\ref{5bcontrol}).

   We will consider all possible continuous matrices $A(x)$ whose elements
are  polynomial functions of x for which the system
(\ref{5control}) is equivalent to the system (\ref{1control}) . We
will denote the set of all such matrices by $\cal A$.

   Let us solve the problem \ref{1probcalopt}.

\begin{thm} In order that a domain D consisting of the
whole trajectories of the system (\ref{1control}) i.e.
$x(\cdot,x_0,t_0) \in D, x_0 \in D,$ for all $t>t_0$ is
a region of asymptotic stability it is necessary and sufficient that there is
a matrix $A(\cdot) \in \cal A$ of the system
(\ref{5control}) in the domain D whose eigen-values have negative
real parts at any point $x \in D, x \neq 0.$
\label{1thmcontrol}
\end{thm}
{\bf Proof. Necessity}. Let a domain D be asymptotic stable.
Consider any trajectory $x(\cdot,x_0,t_0), x(t_0)=x_0.$ There is
such a transformation $\xi=X(x)$ of $\Re^n$ in a neighborhood of a
point $x_1=x(t_1)$ that transforms the system (\ref{1control}) to
the differential system
\be
                        \stackrel{\cdot}{\xi}=B(\xi)\xi
\label{6control} \ee where the matrix \(B\) has the eigen-values
with negative real parts.

Indeed for any point $x(t_1,x_0,t_0), t_1>t_0$ there is a linear
transformation \( T_x \) defined in a small neighborhood of the
point \(x(t_1,x_0,t_0) \) and that is the linear part of the
transformation $X(.)$ in this neighborhood, so that the system
(\ref{1control}) is transformed by this transformation to the
differential system (\ref{6control}) and the vector $$
                        x_1 = \stackrel{\cdot}{x}(t_1,x_0,t_0)
$$
is transformed to the vector
$$
                      \xi_1 = \stackrel{\cdot}{\xi}(t_1,x_0,x_0) .
$$ Any differential equation can be defined by its current of
tangent vectors. If we choose an asymptotic stable linear
differential system like the system (\ref{3control}) which current
of tangent vectors is close to current of tangent vectors of our
system locally for each point $x$ of the domain D then the
eigen-values of the matrix of this asymptotic stable linear
differential system will have negative real parts. Consequently,
the above mentioned transformation $X(x)$ exists.

   As soon as the system (\ref{6control}) is transformed to the system
$$
                       \stackrel{\cdot}{x} =T_x^{-1} B(\xi_1) T_x x
$$ under the linear transformation $$ \xi=T_x \, x      $$ in a
neighborhood of the point $\xi_1$ and all locally linear
transformations of the system (\ref{1control}) have the form
(\ref{5control}), then the matrix \( A(\cdot) \in \cal A \) exists
for which
\be
A(x(t_1, x_0))= T_x^{-1} B(\xi_1) T_x \label{6acontrol} \ee and
the eigen-values of the matrix \(A(x(t_1))\) are equal to the
eigen-values of the matrix \(B\).

Since the transformations $T_x$ and $T_x^{-1}$ are continuous with
respect to $x,$ the matrix $A(.)$ is also continuous with respect
to $x \in D$ .
\\
   {\bf Sufficiency}. Let the system (\ref{1control}) admit such a
transformation to the form (\ref{5control}) that all eigen-values
of a matrix $A(x)$ have negative real parts for all $x \in D, x
\neq 0.$ We prove that the domain D is a region of asymptotic
stability.

  Fist we remark that in general case not any solution of the system
(\ref{1control})
$$        x(t)=x(t,x_0,t_0), \;\; x(t_0)=x_0   $$
can be represented in the form
\be
\stackrel{\cdot}{x}(t)=e^{\int_{t_0}^t A(x(\tau))d\tau} x_0
\label{7control}
\ee
where integration is taken along the integral curve $x(.).$

Instead of (\ref{7control}) we will use the sequential
approximations $\{ x_k(t) \}$       of the system (\ref{1control})
on the segment $[t_0, t]$ having the form
\be
x_k(t)=e^{A(x_{k-1})(t-t_{k-1})}x_{k-1}(t)
\label{7acontrol}
\ee
where $\{t_i \}$ is a subdivision set of the segment $[t_0,t]$ and
\be
  x_{i+1}=e^{A(x_i)(t_{i+1}-t_i)}x_i,\;\;i \in 0:(k-1).
\label{7bcontrol}
\ee
If we substitute (\ref{7bcontrol}) into (\ref{7acontrol}) we will get
\be
x_k(t)=e^{A(x_{k-1})(t-t_{k-1})}e^{A(x_{k-2})
(t_{k-1}-t_{k-2})}...e^{A(x_0)(t_1-t_0)}x_0.
\label{7ccontrol}
\ee
Since $A(x_i) \longrightarrow_{\begin{array}{c} x_i \rar x(\tau)\\
t_i \rar \tau \end{array}}       A(x(t)),$ where $x(.)$ is a solution of
(\ref{5control}), we have
\be
x_k(\tau) \Longrightarrow x(\tau)
\label{7dcontrol}
\ee
uniformly on $\tau \in [t_0,t]$ when $k \rar \infty$ and
$$
      \bmax_{i \in 1:k} \mid t_i-t_{i-1} \mid \rar_k 0.
$$ It is obvious  that all matrices $A(x_i), i \in 0:(k-1),$ have
the eigen-values with negative real-valued parts for enough big k.

  Denote by $\lambda_i(x_j)$ the eigen-values of the matrix
$A(x_j) , i \in 1:n, j \in 0:k-1.$ Then a number $C>0$ exists that
$$ \Vert
e^{A(x_{k-1})(t-t_{k-1})}e^{A(x_{k-2})(t_{k-1}-t_{k-2})}...
e^{A(x_0)(t_1-t_0)}x_0 \Vert \leq C e^{\lambda(t)(t-t_0)}\Vert x_0
\Vert $$ where $$        \lambda(t)= \lim_{k \rar
\infty}\max_{\begin{array}{c} j \in 0:k-1, \\
     i \in 1:n \end{array}} Re \; \lambda_i(x_j).
$$
   Since $ \lambda(t)<0$ for any $t > t_0$,
$$
  \lambda(t)(t-t_0)  < 0.
$$
Let
$$
-a < \lambda(t)(t-t_0)       < 0, a >0.
$$
This relation means that the trajectory $x(\cdot)$ can not go to a stationary
point or be a stationary orbit. Indeed if the trajectory $x(.)$ has a stationary
point or is a stationary orbit then
$$
  \lim_{t \rar \infty}       \lambda(t)(t-t_0) = - \infty.
$$
 Consequently,
$$                     \pr x(t, x_0, t_0) \pr \lrar 0                  $$
when $t \lrar \infty.$ We obtain the contradiction.
Thus, the origin of coordinates can be only a stationary point.
The sufficiency and the theorem are proved. $\Box$

\begin{rem}
It follows from the Necessity of the Theorem  \ref{1thmcontrol}
that there is one-to- one correspondence between the matrices
$A(x) \in \cal A$ of the system (\ref{5control}) and locally
linear transformations $T_x$ of the system (\ref{1control})
defined in a neighborhood of the point $x$(see (\ref{6acontrol}))
. Indeed, there is a matrix $A(.)$ for any transformation $T_x$
and, on the contrary, there is own coordinate system for any
matrix $A(x) \in \cal{A},$ in which the equation (\ref{6acontrol})
is true.
\end{rem}

\begin{cor}.  It follows from the Necessity of the Theorem  \ref{1thmcontrol} that
the degree on x of any element of a matrix $A(x)$ which we denote
by $a_{ij}(x),j \in 1:n,$ that is a polynomial of x, does not
exceed the degree of the vector-function $f (\cdot).$
\end{cor}
   {\bf Proof} The sum of the coefficients $a_{ij}(x) \cdot x_j$ over j is
equal to zero for those elements that do not belong to
$f_i(\cdot).$ And on the contrary the sum of the coefficients
$\al_{j,l_1,l_2,...,l_n}^{(i)}$ over $j$ is not zero for those
elements that are the terms of $f_i(\cdot).$ It follows from here
that choosing $\al_{j,l_1,l_2,...,l_n}^{(i)}$ from the system
(\ref{5acontrol}) we can only consider the elements of the vector-
polynomial $f(\cdot). \,\,\, \Box
$

\begin{cor}. For asymptotic stability of a solution $x(.,x_0,t_0),\,\,x(t_0)=x_0, $
of the differential equation (\ref{1control}) it is sufficient
that there was $\del(t_0,\eps)>0$ for any $t_0,\eps>0$ that for
$\pr x_0 \pr < \del$
 $$ \pr x(t,x_0,t_0) \pr
\lrar 0 $$ when $t \lrar \infty. $ \\
\end{cor}
   {\bf Proof}. As soon as the eigen-values of a matrix $A(\cdot)$
have negative real parts in any neighborhood of the origin of
coordinates from which any solution tends to zero-vector , then
for any k we have from (\ref{7ccontrol}) that $$ \pr x_k(t) \pr
\leq \pr x_0 \pr < \eps  $$ and in limit on k $$ \pr x(t,x_0,t_0)
\pr \leq \pr x_0 \pr < \eps $$ i.e. all solutions are in an $\eps$
-neighborhood of the origin of coordinates. The latter means
stability. The Corollary is proved. $\Box$

\vspace{1.5cm}

\section{To problem about center and focus}

\vspace{1cm}

   The following problem is interesting from technical point of view:
to recognize focus or center i.e. to give the conditions when
trajectories turn around and go to the origin of coordinates or
remain closed.

From the start we will give the definitions of focus or center.
Beforehand we define trajectories turning around some ray with the
initial point $0\in \Re^n$ infinitely often.

\begin{defi}
Let us say that a trajectory $x(t,x_0,t_0)$ turns  around a ray
$l \in \Re^n$ with the initial point $0\in\Re^n$ infinitely often
if the radius vector $r(t)=x(t,x_0,t_0)$ forms with the ray $l$
an angle $\varphi(t)$  taking all values from a segment
$[\varphi_1,\varphi_2]$ infinite number of times.
\end{defi}

\begin{defi}
We will call the point $0 \in \Re^n$  focus of the system
(\ref{1control}) if

1) the point $0$ is asymptotic stable ;

2) any solution $x(t,x_0,t_0)$ of the system $(\ref{1control})$
   turns around a ray $l\in\Re^n$ with the initial point $0\in\Re^n$
   infinite number of times.

\end{defi}

\begin{defi}
We will call the point 0 a center of the system (\ref{1control})
if

1) the point 0 is stable ;

2) all solutions $x(t,x_0,t_0)$ of the system $(\ref{1control})$
   remain in some neighborhood of the point 0 and are closed loops.
\end{defi}

\begin{defi}
We will call the point $0 \in \Re^n$ a node of the system (\ref{1control})
if

1) the point $0$ is asymptotic stable ;

2) there are not a ray $l\in\Re^n$ with the initial point $0\in\Re^n$
   and a solution $x(t,x_0,t_0)$ of the system $(\ref{1control})$
   turning around the ray $l$ infinite number of times.

\end{defi}

   \begin{thm} In order that the point $0 \in \Re^n$  is a focus of the
system (\ref{1control}) with a right polynomial part $f(.), f(x) \neq 0$
for $x \neq 0$ it is necessary
and sufficient that
\begin{enumerate}
\item There is such a continuous matrix $A(\cdot) \in \cal A$ of the
system (\ref{5control}) that all its eigen-values at any point $x,
x \neq 0,$ from some neighborhood $D$ of the origin of
coordinates, consisting from the whole trajectories, have negative
real parts and non-zero imaginary parts.
\item There is no matrix
$A(\cdot) \in \cal A$ of the system (\ref{5control}) with the
negative real-valued eigen-values at all points $x \in D, x \neq
0.$
\end{enumerate}
\label{2thmcontrol}
\end{thm}

   {\bf Proof. Necessity.} Let the point $0$ be a focus for the system
(\ref{1control}). There is such a transformation $\xi=X(x)$ so
that the system (\ref{1control}) can be rewritten in a
neighborhood of some point $x_1, x_1 \neq 0,$ in the form $$
\stackrel{\cdot}{\xi} = B(\xi)\xi       $$ where $B(\xi)$ is a
matrix whose eigen-values have negative real parts and nonzero
imaginary parts.

    Consider any trajectory $x(\cdot, x_0,t_0), x(t_0)=x_0.$ Then
for any point $x_1 \in x(t,x_0,t_0), x_1 \neq 0,$ there is a
transformation $T_{x_1}$ in a neighborhood of the point $x_1$ that
the system (\ref{1control}) can be rewritten in the form indicated
above.

If the trajectories of the system (\ref{1control})  do not turn
around any ray $l\in\Re^n$ with the initial point $0\in\Re^n$
infinite number of times then such system can be transformed by
some continuous transformation that corresponds to some matrix
$A(\cdot) \in \cal A$ that is the matrix of the system
(\ref{5control}) whose eigen-values are negative real-valued
numbers at all points $x \in D$ (theorem \ref{1thmcontrol}). The
latter is impossible according to the Condition 2 of the theorem.

   {\bf Sufficiency.} Let the conditions 1,2 of the theorem be true. Prove
that the point $0$ is a focus of the system (\ref{1control}).

   Any solution $x(\cdot, x_0,t_0), x(t_0)=x_0,$ is the limit of some sequence
$x_k(t)$ on $k \rar \infty$ obtained from $(\ref{7ccontrol}).$
According to the conditions all matrices $A(x_j), j \in 1:k,$ have
the eigen-values $\lambda_l(x_j), j \in 1:k, l \in 1:n$ with
negative real-valued parts $a_l(\cdot)$ and nonzero imaginary
parts $b_l(\cdot): \lambda_l(x)=a_l(x)+i b_l(x), i^2=-1,\,\,\, l
\in 1:n,$ that are nonzero at any point $x \in D.$ It follows from
here that the point $0$ is either focus or node.

   We will prove that under the Condition 2 all trajectories turn around
a ray $l\in\Re^n$ with the initial point $0\in\Re^n$ infinite
number of times. If it is not true, then there is such a
transformation which corresponds to a matrix $A(\cdot) \in \cal A$
whose eigen-values are negative real-valued numbers for $x \in D,
x \neq 0.$ This contradicts the Condition 2. The sufficiency and
the theorem are proved. $\Box$

\begin{cor} If instead of the Condition 2 of the Theorem \ref{2thmcontrol}
we require that there is a matrix $A(.) \in \cal A$ of the system
(\ref{5control}) with the eigen-values $\lambda_i(x), \,\,\,i \in
1:n,$ described in the condition 1 for which $$        Im \;
\lambda_i(x)>a, \,\,\, a>0     $$ or $$        Im \;
\lambda_i(x)<a, \,\,\, a<0,     $$ where $Im \, \lambda_k(.)$
denotes the imaginary part of $\lambda_k(.),$ for $x \in D, x \neq
0,$ then all trajectories turn around a ray $l\in\Re^n$ with the
initial point $0\in\Re^n$ infinite number of times. \label{1cor}
\end{cor}

{\bf Proof.} From the representation of any solution of the system
(\ref{1control}) as a limit of sequence $\{x_j(.)\}$ from
(\ref{7ccontrol}) it follows that a solution $x(\cdot, x_0,t_0)$
has a finite number of turns around any ray $l\in\Re^n$ with the
initial point $0\in\Re^n$ if and only if for any $i \in 1:n$ the
sum $$ \sum_{j \in 1:k} Im  \lambda_i(x_{j-1})(t_j-t_{j-1}) , $$
where $x_j(.)$ was defined in (\ref{7ccontrol}), or in the limit
on $k \rar \infty$ the integral $$
  \int_{t_0}^{\infty} Im \; \lambda_i(x(\tau))d\tau
$$ is convergent. But under the condition of Corollary \ref{1cor}
this integral is divergent. The corollary is proved. $\Box$

   Mathematicians were concerned about the problem of recognition of center
ore focus for a long time. The theorem given below for two
dimensional system with a right polynomial part has the point
$0=(0,0)$  as a center.

\begin{thm} In order that the point $0=(0,0)$ be a center for the
two-dimensional system (\ref{1control}) it is necessary and
sufficient that there was a matrix $A(.) \in \cal A$ of the system
(\ref{5control}) whose eigen-values are non-zero imaginary numbers
in a neighborhood $S$ of the origin of coordinates, consisting
from whole trajectories, where $f(x) \neq 0$ for $x \neq 0, x \in
S.$ \label{thm32}
\end{thm}
   {\bf Proof. Necessity.} Let the point $0=(0,0)$ be a center of the system
(\ref{1control}). Prove that the condition of the theorem is true.

   Take such a small neighborhood $S$ of the origin of coordinates where
all trajectories starting in $S$ are closed loops. Consider any
trajectory $x(\cdot,x_0,t_0), x(t_0)=x_0.$ Take a point
$$x_1=x(t_1,x_0,t_0) .$$
   There is such a transformation $T_x$ in a small enough neighborhood of the
point $x_1, x_1 \neq 0,$ that the system (\ref{1control}) can be
rewritten in the form $$          \stackrel{\cdot}{\xi} =
B(\xi)\xi    $$ and the matrix $B(.)$ has the imaginary
eigen-values for all $\xi \in S.$        As was proved above
(Theorem \ref{1thmcontrol}, the Proof of the Necessity) the matrix
$B(\cdot)$ corresponds to some matrix $A(\cdot)$ that like the
matrix $B(\cdot)$ has the imaginary eigen-values. The necessity is
proved.

   {\bf Sufficiency.} Let the condition of the theorem be true. Prove that the
point $0=(0,0)$ is the center.

   As soon as any transformation $T_x$ is defined in a small neighborhood
of the point $x$, then the above statement is not sufficient for
the point $0=(0,0)$ to be a center. The point $0=(0,0)$ may happen
to be a focus convergent (nonconvergent) with a finite number of
turns around any ray $l\in\Re^n$ with the initial point $0$. But
it is impossible because from the limit of the sequence
(\ref{7ccontrol}) we conclude that for any k $$ \pr x_k(t) \pr =
\pr x_0 \pr. $$ From (\ref{7dcontrol}) it follows that the
solution $x(\cdot,x_0,t_0)$ can't go to the origin of coordinates.

A set of closed and not closed loops could occur. From the
representation of an integral curve $x(\cdot, x_0,t_0)$ in the
form (\ref{7ccontrol}) and (\ref{7dcontrol}) it follows that there
is not a matrix $A(\cdot) \in \cal A$ with the imaginary
eigen-values at all points of a small neighborhood of the origin
of coordinates. The latter contradicts the condition of the
theorem.

   There is another case when for all $i \in 1:n$
the sum
$$
\sum_{j \in 1:k} Im  \lambda_i(x_{j-1})(t_j-t_{j-1})
$$
or in the limit on $k \rar \infty$ the integral
$$
  \int_{t_0}^{\infty} Im \; \lambda_i(x(\tau))d\tau
$$
has a finite value where $\lambda_i(\cdot),\,\,\,i \in 1:n$, are
the eigen-values of the matrix $A(\cdot).$ In this case the
trajectory $x(\cdot,x_0,t_0)$ goes to some stationary point $\hat{x}
\in S, \hat{x} \neq 0.$  It means that $f(\hat{x}) = 0.$ This
contradicts the condition of the theorem.       Consequently, the
point $0$ is the center. The sufficiency and the theorem are
proved. $\Box$

\begin{cor} If a matrix $A(.) \in \cal A$ of the system (\ref{5control})
exists with the imaginary eigen-values $\lambda_i(x)$ for all $x
\in S$  where $S$ is a neighborhood of the origin of coordinates
$0$, consisting from whole trajectories, and $$             Im
\,\, \lambda_i(x) > 0          $$ or $$             Im \,\,
\lambda_i(x) < 0         $$ for all $x \in S, x \neq 0,$ and $i
\in 1:n$ then the point $0$ is the center of the n-dimensional
system (\ref{1control}). \label{2cor}
\end{cor}
{\bf Proof.} Indeed, if the conditions of the Corollary \ref{2cor}
hold, then the conditions of the Theorem \ref{thm32} hold as well.
$\Box$

   The Corollary \ref{2cor} does not require that $f(x) \neq 0$ for all
$x \neq 0$ from some neighborhood of the origin of coordinates.

\begin{rem} The theorem can be proved for any n-even-dimensional spaces.
It is not difficult to see that for n dimensional spaces with n
odd, $n=2k+1, k$ is a natural number, the system (\ref{1control})
can not have the origin of coordinates as a center. Indeed, there
is no imaginary number among the eigen-values of any matrix
$A(\cdot)$ for odd n.
\end{rem}

\begin{rem} Let us consider the differential system
\be
                \stackrel {\cdot}{x} = f(x,t) \qquad \qquad \qquad
\label{8control}
\ee
where \( x=(x_1,x_2,...,x_n) \in \Re^n \) and
$$
              f(x,t)=(f_1(x,t),f_2(x,t),..., f_n(x,t))^{*}.
$$ The vector-function $f_p(.)$ is a polynomial of $x$ i.e. $$
f_p(x) = \sum_{l_1,l_2,...,l_n \in I_p}
a_{l_1,l_2,...,l_n}^{(p)}(t) x_1^{l_1}x_2^{l_2}...x_n^{l_n} , $$
\( p \in 1:n, l_1,l_2, ... l_n \) are non-negative integers,
$a^{(p)}_{l_1,l_2,...,l_n}(t)$ are continuous real-valued
functions and \(* \) is the transposition sign, $I_p$ is a finite
set of indexes of the polynomial $f_p(.).$ We will assume that
$f(x,t) \neq 0$ for all $x \neq 0$ in some neighborhood of the
origin of coordinates $0$ and $t > t_0.$

In this case  we will denote by $\lambda_i(x,t), i \in 1:n,$ the
eigen-values of a matrix $A(.).$ If we demand that all cited above
statements about the eigen-values $\lambda_i(x,t), i \in 1:n,$ are
true for all $t > t_0$ then we obtain the sufficient conditions
for the cited above theorems and corollaries.
\end{rem}

\begin{prob} For the system
\be
             \stackrel{\cdot}{x} = f(x,t)
\label{9control} \ee where $f(.,.)$ is a polynomial of $x$ and $t$
it is required to find some conditions when the system
(\ref{9control}) is asymptotic stable in a domain $D, 0 \in int
D.$ \label{3probcalopt}
\end{prob}

   The ideas stated above do not apply to the system (\ref{9control})
because the terms of $f(.,.)$ can be unbounded along some solution
$x(\cdot,x_0,t_0)$ in D. If we require
the terms of $f(.,.)$       to be bounded along the trajectories of the system
(\ref{9control}) and that $\pr \stackrel{\cdot}{x}(t,x_0,t_0) \pr$ goes
to zero uniformly on $x_0$ when $t \lrar \infty$, then this system can be
transformed to an equivalent stationary system. This idea will be developed
in following articles.

\vspace{1cm}

\section{Stabilizing control}

\vspace{1.5cm}

   Let us go to the question of finding a stabilizing control.

   Let us consider the differential system
\be
               \stackrel{\cdot}{x} = f(x,u)
\label{10control}
\ee
where $x=(x_1,x_2,...,x_n) \in \Re^n$ is a phase vector,
$u=(u_1,u_2,...,u_r) \in \Re^r$
is a control, $f(x,u)=(f_1(x,u),f_2(x,u),...,f_n(x,u))^*$ is a
vector-polinimial of x and u with constant real-valued coefficients, i.e.
$$
f_p(x,u) = \sum_{l_1,l_2,...,l_n,m_1,m_2,...,m_r \in I_p}
a_{l_1,l_2,...,l_n,m_1,m_2, ...,m_r}^{(p)}
x_1^{l_1}x_2^{l_2}...x_n^{l_n}u_1^{m_1}u_2^{m_2}... u_r^{m_r},
$$
 \( p \in 1:n,\,\,\,\,l_1,l_2,...,l_n, m_1,m_2, ..., m_r \) are non-negative
integers and $a^{(p)}_{l_1,l_2,...,l_n,m_1,m_2,...,m_r}$ are real-valued
numbers, $I_p$ is a finite set of indexes of the polynomial $f_p(.).$

Let us assume that the zero-vector $0=(0,0,...,0)
\in \Re^n$ is a solution of the system (\ref{10control}) for $u=0 \in \Re^r.$

   \begin{defi} The control $u(x)=(u_1(x),u_2(x),
...,u_r(x)) \in \Re^r$ is called stabilizing in a domain $D \in
\Re^n, \,\, 0 \in intD,$ for the system (\ref{10control}) if any
solution of (\ref{10control}) $x(t)=x(t,x_0,t_0,u(x(t)))$
satisfies the limit
\be
           \pr x(t,x_0,t_0,u) \pr \lrar 0
\label{11control}
\ee
when $t \lrar \infty, \,\,\, x(t_0)=x_0 \in D.$
\end{defi}
   As stated above the condition (\ref{11control}) is sufficient for the
system (\ref{10control}) to be asymptotic stable in D, i.e. the
zero solution of the system (\ref{10control}) is asymptotic stable
and the limit (\ref{11control}) is true for any initial point
$x(t_0)=x_0 \in D.$

\begin{prob} It is required to find such a stabilizing control $u=u(x)$
in a given domain $D,0 \in int D, x \in D$ that for any solution $x(t,x_0,t_0,u),
x(t_0)=x_0$ of the system (\ref{10control}) the correlation (\ref{11control})
was true.
\label{4probcalopt}
\end{prob}

   For linear system
\be
            \stackrel{\cdot}{x} = Ax+Bu,
\label{12control}
\ee
where $B[n \times r]$ is a matrix of amplification coefficients, a stabilizing
control $u \in \Re^r$ can be chosen in the form $u=Cx$ so that
(\ref{11control}) is true.

\begin{thm} \cite{zub1},\cite{kras1}. If the rank of the system
$$
     B,\,\,AB,\,\,A^2B,\,\,...,\,\,A^{n-1}B
$$
is equal to n then we can always construct a stabilizing control in $\Re^n$
in the form
$$                    u=Cx,            $$
where $C[r \times n]$ is some matrix.
\end{thm}
   Since the matrix $B$ is defined by technical essence of the system,
we can choose it constructing a system by ourselves . Therefore the
condition of the theorem can be always satisfied.

The theorem was written for analogy and comparison with the
results below.

   Our goal is to construct a stabilizing control $u(x)$ that a domain $D$
were a domain of asymptotic stability. In common case it is not always can be
done. But if we change the system (\ref{10control}) a little bit the problem
can be solved. Instead of the equation (\ref{10control})
we will consider the equation
\be
     \stackrel{\cdot}{x}=f(x,u)+\varphi(u)
\label{13control} \ee where $\varphi(\cdot)$ is a
vector-polynomial $\varphi(u)=(\varphi_1(u),
\varphi_2(u),...,\varphi_n(u))^{*}$ with the degree not bigger
than the degree of the function $f(z)$ as a function of
$z=(x,u),$ $$ \varphi_p(u) = \sum_ {i_1,i_2,...,i_r \in M_p}
\,\,b_{i_1,i_2,...,i_r}^{(p)} u_1^{i_1}u_2^{i_2}... u_r^{i_r} $$
where $p \in 1:n$ and $b_{i_1,i_2,...,i_r}^{(p)}$ are constant
real-valued numbers, $i_1,i_2,...,i_r$ are non-negative integers,
$M_p$ is a finite set of indexes of the polynomial $\varphi_p(.).$

   In practice it is possible to construct $\varphi(\cdot)$ because we choose
a control $u(\cdot)$ ourselves.

  We will look for a stabilizing control $u=u(x)$ in form of polynomial in
 $x.$
\begin{thm} For any domain $D, 0 \in intD$, the vector-polynomials $u(.)$ and
$\varphi(.)$ can be chosen such that \\ 1) $D$ is a region of
asymptotic stability for the differential system
(\ref{13control}); \\ 2) the degree of $u(x)$ does not exceed the
degree of the vector-polynomial $f(x,u)$ as a function of x; \\ 3)
the degree of $\varphi(.)$ is not bigger than the degree $f(.)$ as
a function of \\ $z=(x,u).$ \label{4.2controlthm}
\end{thm}
   {\bf Proof.} We will prove that a vector polynomial $u(x)$ can       always
be chosen satisfying the following conditions \\
1) the degree of $u(x)$ does not exceed the degree of the vector-polynomial
$f(x,u)$ as a function of $x$;\\
2) $u(x)$ is a stabilizing control for (\ref{13control}) in the domain $D.$

   We will use the results obtained before (Theorem \ref{1thmcontrol}).

   Let us substitute in (\ref{13control}) the control $u=u(x)$ in the form
of a vector-function of $x.$ The system (\ref{13control}) is rewritten as
\be
\stackrel{\cdot}{x} = \hat{f}(x)
\label{14control}
\ee
or
$$
\stackrel{\cdot}{x} = \hat{A}(x)x.
$$
   We can write conditions for the matrix $\hat{A}(x)$ to have the eigen-values
with negative real parts.

   There is no difficulty calculating the degree of the function $f(x,u)$
and the function $\varphi(u)$ as a function of $x$ after substituting $u=u(x).$

   Let us denote by $l_x$ and $l_u$ the degrees of the function $f(x,u)$ in the
variables $x$ and $u$ correspondently. Then after substituting $u=u(x)$
the degree of the function $f(.,.)$
as a function of $x$ is not bigger than $l_x +l_u l_x$ and the degree of the
function $\varphi(.)$
as a function of $x$ is not bigger than $l_x(l_x + l_u).$ It is easy to see that
$$        l_x+l_x l_u \leq l_x(l_x+l_u).   $$

   Having chosen the function $\varphi(\cdot)$,
the coefficients of the vector polynomial $\hat{f}(x)$ can be
chosen so that the matrix $\hat{A}(x)$  has the eigen-values with
negative real parts. The theorem is proved. $\Box$

   In an analogous way it can be proved the following theorem.
\begin{thm} For the differential system (\ref{13control}) with even n the
vector-polynomials $u(\cdot)$ and $\varphi(\cdot)$ can be chosen
such that \\ 1) the origin of coordinates 0 is the center for the
system (\ref{13control});\\ 2) the degree of $u(x)$ does not
exceed the degree  $f(x,u)$ as a function of $x$\\ 3) the degree
of $\varphi(.)$ is not bigger than the degree of $f(.)$ as a
function of $z=(x,u).$ \label{4.3controlthm}
\end{thm}

   This theorem can be used in physics of plasma specially for stabilization
of plasma in reactor.

\begin{rem}
In the case when the coefficients $a^{(p)}_{l_1,l_2,...,l_n,m_1,m_2,...,m_r}$
depend on t the coefficients $b_{i_1,i_2,...,i_r}^{(p)}$ will
depend on t as well. We can try to find these coefficients that the conditions
of the theorems of the previous section were true for all $t > t_0.$
\end{rem}

\vspace{1cm}

\section{One aspect of application}

\vspace{1.5cm}

   We will study for instance how to choose a control $u(\cdot)$ so that a
given domain $D$ consisting from the whole trajectories of a differential system
were a region of asymptotic stability. For that we have to solve the following
optimization problem.

Let us substitute a vector-control $u(\cdot)$ in a form of
polynomial of $x$ with a degree $m$ not bigger than the degree of
the the vector-polynomial $f(x,u)$ as a function of $x$ in the
equation (\ref{13control}) and rewrite (\ref{13control}) in the
form (\ref{14control}). Denote the eigen-values of a matrix
$\hat{A}(\cdot)$ in (\ref{14control}) by
$\lambda_j(x,\al_{ij},d^{(l)}_{i_1,i_2,...,i_n},
b_{i_1,i_2,...,i_r}^{(p)})$ where $\al_{ij}$ and
$d^{(l)}_{i_1,i_2,...i_n}$ are the coefficients of the matrix
$\hat{A}(\cdot)$ and the vector-function $u(\cdot)$
correspondingly, $b_{i_1,i_2,...,i_r}^{(p)}, \, p \in 1:n,$ are
the coefficients of the vector-function $\varphi(.)$. Then our
problem is reduced to the following optimization problem: to find
out such continuous functions $\al_{ij}(x),i,j \in 1:n, \, x \in
D$ and the numbers $d^{(l)}_{i_1,i_2,..,i_n}, \,
b_{i_1,i_2,...,i_r}^{(p)},\, p \in 1:n,$ that the following
correlation
\be
Re \, \lambda_j(x,\al_{ij}(x),d^{(l)}_{i_1,i_2,...,i_n}, b_{i_1,i_2,...,i_r}^{(p)}) < 0
\label{15control}
\ee
was true for all $x \in D, x \neq 0,$ where $\al_{ij}(x), i,j \in 1:n,$ are
connected with
each other by linear equations for each $x \in D,x \neq 0,$ by other words we
should solve the next problem
$$
     \inf_{\begin{array}{c} d^{(l)}_{i_1,i_2,..,i_n},l \in 1:r \\
      {b_{i_1,i_2,...,i_r}^{(p)}, \, p \in 1:n } \end{array} } \,\,\,
     \sup_{\begin{array}{c}
      {j \in 1:n} \\ {x \in D \setminus B^n_{\del}(0) }  \end{array}} \,\,\,
      \inf_{c_{ij}} \,
      Re \, \lambda_j(x,\al_{ij},d^{(l)}_{i_1,i_2,...,i_n},
      b_{i_1,i_2,...,i_r}^{(p)}) < 0  \,\,\, \forall \, \del >0 ,
$$
where $B^n_{\del}(0)=\{z \in \Re^n \mid \pr z \pr \leq \del \}, \del -       $
is any sufficient small number for which $D \supset B^n_{\del}(0).$

It is obvious that the inequalities (\ref{15control}) can be
replaced by an equivalent system of inequalities for coefficients
of the characteristic polynomial of the matrix $\hat{A}(.).$ To
solve this system  it is easier than to find the eigen-values of
the marix $\hat{A}(.)$.

\end{document}